# On the history of analysis.
# The formation of definitions

G.Sinkevich

Saint Petersburg State University of Architecture and Civil Engineering,
Vtoraja Krasnoarmejskaja ul. 4, St. Petersburg, 190005, Russia
galina.sinkevich@gmail.com

**Absract**. The mathematical analysis was conceived in XVII century in Newton and Leibniz works. The problem of logical rigor in definitions was considered by Arnauld and Nicole in "Logique ou l'art de penser". They were the first, who distinguished between the volume of the concept and its structure. They created a tradition which was strong in math till XIX century, especially in France. The definitions were in binomial nomenclature mostly, but another type of definition appears in Cantor theory – it was the descriptive definition. As it used to be in human sciences, first the object had only one characteristic, then as the research continued it enriched with new characteristics so we have full-fledged conception. In this way the mathematics acquired its own creativity. In 1915 Lusin laid down a new principle of the descriptive theory: a structural characteristic is done, the analytical form had to be found. New schools of descriptive set theory appeared in Moscow in the first half of the 20th century.

Mathematical analysis took its rise in Newton's and Leibniz' works. For Newton, existence of a derivative was evident since a speed of movement has to exist, his concept of limit containing no restrictive conditions. For Leibniz, everything that is possible tends to exist .

The problem of logical rigor of definitions arose with transfer from Latin to national languages. Leibniz was dreaming of creation of a universal language which would make it possible to abrogate Latin having replaced it with characteristic (a system of signs of thoughts). In the French language, this problem was first laid down in the 17th century in works of philologists, philosophers and mathematicians who were members of a scientific circle at Port-Royal-des-Champs Convent. Blaise Pascal attended this circle.

Port-Royal Convent was founded in the 12th century not far from Paris (in the 17th century it was 6 leagues away). The vicinity of Paris determined the aristocratic population of the convent – daughters of noble families were sent to this place. Customs in this convent were not puritan, relatives of novitiates and nuns used to stay in its guest house for long periods. There was a good school at the convent, teachers would write schoolbooks and print them at the printing house of the convent. In 1602, Angélique (Jacqueline) Arnauld became their abbess. Her brother, Antoine Arnauld, prominent theologian and mathematician, theologian Pierre Nicole, and philologist Claude Lancelot were among teachers of the school. Blaise Pascal often visited his cloisterer sister Jacqueline Pascal. He later settled here in the guest house. Aristocrats, military officers, writers and scholars living there permanently or coming from Paris formed a philosophers circle around them. In 1660, Arnauld and Lancelot published a Grammar Book of Port-Royal (Grammaire générale et raisonné de Port-Royal) which set forth a new language conception.

Latin as the only possible language of science gradually yielded its position. The authors of the Grammar Book were discussing the issue of a common logical basis for all languages which corresponds to the structure of the thought. They also discussed the degree of accuracy of conveyance of the meaning when translated from one language to another. In other words, when translating a text into various languages, it is only the logical part that we keep unchanged. For example, this article was translated from Russian, the meaning has been conveyed, however, something imperceptible was lost.

The authors of "Logic" aspired to identify the fundamental structures of human conscience. In "The Logic or Art to Think" (Logique ou l'art de penser) Arnauld and Nicole look into this matter in terms of notions in science at large and in mathematics in particular. The text incorporates critically revised speculations of Euclid, Stevin, Descartes, Pascal, and many others.

"Logique ou l'art de penser" implied liberation from verbal forms, search of the initial meaning, which is the common component of all languages. It was for the first time that they differentiated the volume of the concept and the structure thereof, definition of the idea (definition nominis) and definition of the real thing (definition rei). This method had to become not only the method of substantiation, but the method of research as well. Mathematicians were blamed to lack rigor, to use the grounds anticipation principle (petitio principii). Arnauld and Nicole write in Section "The Fifth Drawback – Not to Think of the Natural Order: This is the greatest drawback of geometricians. They believe they need not adhere to any rule, except that first provisions have to manifest the previous ones. Therefore, neglecting the rules of the right method which implies that one should always start with the simplest and most general and thereafter pass over to more complex and more particular, they speak indiscriminately of lines and areas, triangles and squares, prove properties of simple lines through figures, and violate numerous other rules [natural in their nature], which corrupt this wonderful science. In Euclid's "The Bases", one would find this drawback everywhere".

The principle of invariance of the volume and structure of the notion was asserted in the course of the discussions. For example, when the term "even number" is used, which means a number which can be divided into 2, we can at any time in the course of the discussion similarly use the words "a number which can be divided into 2" instead of the words "even number" [Арно, Николь, с. 337–338].

Arnauld and Nicole compare the two definitions of the number provided by Euclid and Stevin. According to Stevin, a number is what expresses an amount of anything. Therefore, a unit is a number. According to Euclid, a number is a set composed of units. Therefore, a unit is not a number. Arnauld and Nicole explain these contradictions by the fact that the definition of the word which cannot be disputed is substituted by the definition of the thing, which can be often rightfully disputed .

Mathematicians addressed these issues again in the 19th century when the problem of rigor in mathematics arose. French mathematicians retained the Port-Royal's tradition of logic.

A.Cauchy (1879–1857) started regulating analysis. He was also the first to try and prove the existence of mathematical objects, e.g., partial solution of a differential equation without constructing them. Cauchy aspired to make definitions as rigorous as required in geometry without recourse to evidence proceeding from generality of algebra.

Euler gave the concept of function, which was used in the 18th century. It meant that the function must be defined by a single formula. In the 19th century it was replaced by the definition of the Dirichlet-Lobachevsky. One-to-one correspondence involves both a necessary and sufficient the totality of properties of functions hence it follows the analytical form.

In the middle of 19th century mathematics was in need of the theory of real numbers. Fourier series appeared, class of discontinuous and non-differentiable functions was expanded. This gave rise to the need to evaluate the set of discontinuity points and the ability to disregard such sets. By the middle of the XIX century already there are new ideas about the need to analyze sets of points, but there are no the words to express them yet. For example, in works of Riemann we see that individual points and curves can be present only in terms of the final sense, and two "individual curves" cannot intersect except at a finite number of points. Riemann with virtuoso dexterity avoids talking about point sets, "in order to avoid a lengthy review of the details, which are not essentially necessary"[Риман, с. 27–28].

It was a new nature of teaching, the main part of mathematical foundations was transferred in the lectures. Lections referred not to intuitive analogies in physical and geometric images, but to the logical perception mostly. This leads to the arithmetization of analysis.

In the 1830s, Bernard Bolzano (1781-1848). attempted to create a theory of real numbers. His theory contained a concept of section together with the concept of the least upper bound, but also contained and actually infinitely large and infinitely small numbers. His concept of least upper bound was not built entirely - Bolzano only shows that the presence of such a boundary does not involve to a contradiction.

In 1821 Cauchy in his "Algebraic analysis " defined the irrational number as the limit of a convergent sequence, but he introduced neither ordering relationship, nor rules of procedure for irrational numbers.

In 1869, the French mathematician Charles Measure (Meray, 1835-1911) had published his theory of real numbers. It was fully expounded in his "Analysis" in 1872 [Méray]. This construction was based on the concept of numerical sequences (variants) satisfying the Cauchy criterion (Although Cauchy published this criterion in 1821, Bolzano was the first who formulated it in 1817). Meray introduced equivalence relation for numerical sequences , and every irrational number was the limit of such sequence. If the sequence converges to a point that does not allow a precise definition, Meray called such limit sequence as fictitious. Meray determined the ordering relationship and arithmetic operations for such limits as the same for the corresponding sequences. But Meray did not prove the validity of these properties, but only showing it by example. He confined himself this, and did not build a hierarchy of sequences.

*"We will determine all immeasurable numbers, approximating their values with the help of a* $\delta$*, however small it might seem"* [Meray, p. 2].

Meray defines a continuous function as Cauchy, arguing that these inequalities are so natural, he need not justify their validity, it tested repeatedly. In his opinion, all continuous functions are decomposable into a series of integer powers, and other functions and nonpower series need not be considered.

Meray`s language was not easy for the audience, he introduced many own terms for concepts that are really needed, but defined them only for a narrow class of objects. For example, he introduces the concept olotrope function (Olo tropos, ολο τροπος – all the way), which is similar to a definition of the limit of a uniformly convergent series, (Meray italics) such way:

«Let $f(x, y, …)$ is a sum of entire series, $R_x$, $R_y$,... – their convergence domains (circles),) $R_x^\circ, R_y^\circ,...$ – correspondingly smaller positive values. Let $x', y',...,x, y,...$ are arbitrary variants in this circles. Its centres are in points $O_x, O_y,...$ as centres of domains $R_x^\circ, R_y^\circ,...$, so all differences $x'-x, y'-y,...$ will be infinitely small, and the difference $f(x', y',...) - f(x, y,...)$ will tent to zero.

*Consider the sum of the first N elements of the whole (integer powers) series. Replacing options converge on options that converge in a smaller inner circle of convergence, we obtain the convergence option and their equivalence to each other, N - substituted version of the index variables, infinitely converging to some values.*

That is what will be just that ensure an integer series of the property to determine the value of the function, even immeasurable, of the independent variables. We think the different point of view is unwarranted in this theory.

*If $x', y',...$ respectively tend to these limits x, y, …, located in the inner part of the circle of convergence, then $f(x', y',...)$ tends to $f(x, y,...)$.*

We say that a function $f(x, y,...)$ is olotrope on portions $S_x, S_y,...$ of the auxiliary plane defined by geometrically independent variables when for all the values of the variables *x, y, …,* in these areas $f(x+h, y+k,...)$ is a series of integer powers, and *h, k, …* are such that they remain within the domain of convergence, differing by an amount $\delta_x, \delta_y,...$ not exceeding, not all are equal to zero. We call $\delta_x, \delta_y,...$ *olotrope areas* or *olometres* of a function those sites that are running the named conditions.

Immediately from this definition it follows that *f* (*x, y, …*) is in the form of a series in integer powers of $x-x_0, y-y_0,...$, where $x_0, y_0,...$ represent some numerical values of the variables within the region, such that for neighboring points $x$, $y$, …, the modules their differences remain within olometres. Thus, it turns out that instead of $x-x_0, y-y_0,...$ we can write $h$, $k$, …, $f(x_0+h, y_0+k,...)$.

*The sum of integral power series is olotrope function of the variables inside the circle with the same center as the circle of convergence, but smaller* "[Méray, c. 38].

Complexity of presentation, conservatism, reluctance to go beyond the classical analysis early XIX century, limiting the analysis to entire functions, unnecessarily cumbersome and often unsuccessful terms defining new concepts (given as an example above definition of non-predicative) were inherent Meray. His theory has not found acceptance among his compatriots. Only a century later, this concept of number in France was admitted as "the concept of Meray-Cantor"[Синкевич (2012)].

By the middle of the XIX century the influence of tradition of Cauchy in France and Gauss in Germany has decreased.

At the University of Berlin in 1856 Weierstrass began teaching. He has forbidden to publish his lectures, accusing the French mathematicians in reverse "even assistant professor can lecturing instead professor on their synopses". These lectures were published in the late 19th and in the 20th century thanks to the records of his students [Dugac, Weierstrass].

Weierstrass (1815–1897) delivered his first lecture circuit devoted to immeasurable numbers in 1861/62 academic year. Records of his lectures of 1878 are also available. In summer term of 1886, in response to reproaches of L. Kronecker to the effect of insufficient justifiability of lectures on theory of analytic functions, Weierstrass read additional chapters devoted to foundations of the theory of functions [Weierstrass]. By that time, concepts of a number of Cantor, Heine, and Dedekind already appeared. Weierstrass attempts to critically summarize them and align them with the classical concept of a number as a ratio.

Weierstrass had arithmetized analysis, reducing its dependence on the geometry. It was he who introduced the definition of the language of ε-δ (epsilon-delta) that Cauchy used in 1821 but without revealing the relationship between ε and δ. He had his own concept of real number as an aggregat of the units, based on the uniform convergence of absolutely convergent series. His lectures amaze by weakly formalization, accuracy of mathematical thought and wealth image over the text. [Синкевич (2013)]. Weierstrass uses impredicative definitions, such as "we understand the variable as a value, which is determined by an infinite number of variables that correspond to this definition"[Weierstrass, c.56].

Weierstrass notes incompleteness of a field of rational numbers, gives consideration to the difference between concepts of a number and numerical value. According to Weierstrass, a number is a collection, finite aggregate, e.g. in the form of a decimal notation. A point on a line corresponds to each number, however, it is not obvious that a number corresponds to each point. Unlike his colleagues, he defines a real number as a limit of partial sums of absolutely convergent series, noting the need in arithmetization of the concept of a limit. He introduces order and completeness with respect to arithmetic operations.

Weierstrass created his reasoning of the theory of analytic functions. Conceptions he introduces are not global in their nature – they are necessary for his constructions only. He introduces his own concepts of a continuum and connectivity which differ from those of Cantor; for analytic continuation, he simultaneously builds up a chain of open discs, which is equivalent to Heine covering lemma. Weierstrass defines a number so that it would be sufficient to define continuous changes in arithmetical values in their mutual dependence, "that is to say, an arithmetic expression is calculated in such detail that for any accuracy requirement for any amount $t$ a function may be represented with any approximation. It is always possible to find a mathematical expression for a strictly defined continuous function as well." However, if a function represents series, this does not narrow down, this rather expands, opportunities for study of this function, but the series must have a uniform convergence. "For any value of *x* for which a function has been determined, it can in fact be represented."

Note that in the German language had no such influence of the Latin tradition. It has more synonymic flexibility in such basic mathematical terms, as a number, a numeric value, point, position, border, limit - each translator would agree with me, . German language is more conceptual than the French language, it was enriched by the philosophical tradition of Kant, Hegel and Schopenhauer.

Weierstrass showed mathematics as a field of unsolved problems. He created a galaxy of strong researchers, among whom were Cantor, Schwartz, Mittag-Leffler, Kovalevskaya, Klein and many others.

Simultaneously with the appearance of works Meray (1869, 1872 years) in Germany in 1872 came three articles on the concept of a real number. There were works by R. Dedekind, E. Heine and Cantor.

E. Heine (1821–1881) in 1872 in "Die Elemente der Functionenlehre" had defined irrational numbers as classes of equivalent sequences. "The theory of functions is for the most part developed using elementary fundamental theorems, although insightful researches cast some doubt on certain results, as research results are not always well-reasoned. I can explain it by the fact that, although Mr. Weierstrass' principles are set forth directly in his lectures and indirect verbal communications, and in manuscript copies of his lectures, and are pretty widely spread, they have not been published as worded by the author under the author's control, which hampers the perceptual unity. His statements are based on an incomplete definition of irrational numbers, and the geometric interpretation, where a line is understood as motion, often misleading. Theorems must be founded using the new understanding of real irrational numbers which have been rightfully founded and exist, however little they may differ from rational numbers, and the function has been uniquely determined for each value of the variable, whether it is rational or irrational.

Not that I am publishing this work unhesitatingly long since its first and more significant part *About Numbers* has been finished. Apart from complexity with presentation of such topic, I was hesitant about publishing results of the verbal exchange of ideas which contain former ideas of other people and those of Mr. Weierstrass in the first place, so, all is left to do is to implement these results, which is extremely important not to leave any vague issues in my narrative. I am especially thankful to Mr. Cantor from Halle for the discussion which significantly affected my work, as I borrowed his idea of general numbers which form series.

Let us call a numerical sequence a sequence consisting of numbers $a_1, a_2, a_3, ..., a_n, ...$, when for each arbitrarily small non-vanishing number η such *n* number can be found that $a_n - a_{n+v} < \eta$ can be achieved for all whole positive *v*.

Let us assume that for the structure of (rational) numbers $a_1, a_2, ...$, there is such (rational) number $U$ that $U - a_n$ decreases as *n* grows. In this case, $U$ is the limit of *a*.

We will call general numbers, which in particular cases become rational numbers, as first-order irrational numbers. Like irrational numbers are formed from first-order rational numbers *A*, in the same way, second-order numbers $A'$ can be obtained from limits of irrational numbers, whereupon, third-order irrational numbers $A''$ can be obtained from them, and so on. We will designate $A^{(m)}$ irrational numbers of $m+1$ order". [Heine].

R. Dedekind (1831–1916) in 1872 in his work "Stetigkeit und irrationale Zahlen" reviews properties of equality, order, density of a multitude of rational numbers R (numerical body, a term introduced by Dedekind in appendices to Dirichlet's lectures he published). However, he tries to avoid geometric representations. Having defined the relation "larger" or "less", Dedekind confirms its transitivity; existence of an infinite multitude of other numbers between two numbers; and, for any number, breaking down a multitude of rational numbers into two infinite classes, so that numbers of one of them are less than this number and another one whose numbers are greater than this number; and the number which breaks down the numbers as described above may be assigned either to one class or to the other, in which event it will be either the greatest for the first class or the smallest for the second one.

Further, Dedekind reviews points on a line and sets properties for them same as he has just set for rational numbers, stating that a point on the line corresponds to each rational number.

"However, there are infinitely many points on a line which do not correspond to any rational number, e.g., the size of a diagonal line of a square with a unit side. This implies that the multitude of rational numbers needs to be supplemented arithmetically, so that the range of new numbers could become as complete and continuous as a line. Formerly, the concept of irrational numbers was associated with measurement of extended values, i.e. with geometrical representation. Dedekind tends to introduce a new concept by purely arithmetic means, that is, to define irrational numbers through rational numbers:

If the system of all real numbers is split into two classes, so that each number of the first class is less than each number of the second class, there is one, and only one, number which makes this split.

There are infinitely many sections which cannot be made by a rational number. For example, if $D$ is a square-free integer, there is a whole positive number $\lambda$, so $\lambda^2 < D < (\lambda+1)^2$ . Therefore, it appears that one class has no greatest, and the other class has no smallest number to make a section, which makes the set of rational numbers incomplete or discontinuous. If that's the case and the section cannot be made by a rational number, let us create a new, irrational, number which will create the section. There is one, and only one, rational or irrational number which corresponds to each fixed section. Two numbers are unequal if they correspond to different sections. Relations "larger than" or "less than" may be identified between them.

He defines calculations with real numbers. Herewith, he proves the theorem on continuity of arithmetic operations: "If number $\lambda$ is a result of calculations which involve numbers α, β, γ, …, and of $\lambda$ lies in interval L, one can specify such intervals A, B, C (in which numbers α, β, γ,… lie) that the result of a similar calculation in which, however, numbers α, β, γ,…are substituted with numbers of respective intervals A, B, C,…, will always be a number which lies in interval L."[Dedekind].

Definition of real number given by Dedekind more tending to algebra. It is flawless categorically, but it is only suitable for linearly ordered sets, and uncomfortable functionally, Cantor reproached him for it. But Dedekind`s merit was that he claimed that the number as a concept existing only "in the world of our thoughts," and the freedom for mathematicians to determine the number as it is necessary for its immediate objectives.

G. Cantor (1845–1918) was a pupil of Weierstrass in Berlin, a colleague of Heine in Halle and a friend of Dedekind in correspondence and in summer rest. In determining the number Cantor based on fundamental sequences [Cantor (1872)], this idea, as Heine himself admits, belongs to Cantor.

Cantor constructs a set of numerical values currently known as real numbers, supplementing a set of rational values with irrational numbers using sequences of rational numbers he called fundamental, i.e. sequences that meet the Cauchy criteria. Relations of equality, greater, and less are determined for them.

In the same way, it can be asserted, says Cantor, that a sequence can be in one of the three relations to rational number *a*, which results in $b=a$, $b>a$, $b<a$. Consequently, if *b* is a limit of the sequence, then $b-a_n$ becomes infinitely small with growing *n*. Cantor calls the totality of rational numbers domain *A* and the totality of all numerical *b* values domain *B*. Numeric operations common for rational numbers (addition, subtraction, multiplication, and division, where the factor is non-vanishing) and applied a finite number of times can be extended to domain *A* and *B*. In this event, the domain *A* (that of rational numbers) is obtained from the domain *B* (that of irrational numbers) and together with the latter forms new domain *C*. That is to say, should you set a numerical sequence of numbers $b_1, b_2, ..., b_n, ...$ with numerical values *A* and *B* not all of which belong to domain *A*, if this sequence has such a property that $b_{n+m} - b_n$ becomes infinitely small with growing *n* and any *m*, such sequence is said to have a certain limit *c*. Numerical values *c* form domain *C*. Relations of equality, greater than, less than, and elementary operations are determined as described above. However, even a recognized equality of two values *b* and *b′* from *B* does not imply their equivalence, but only expresses a certain relation between sequences to which they are compared.

Domain *D* is similarly obtained from domain *C* and preceding ranges, and domain *E* is obtained from all above domains, etc.; having completed λ of such transformations, domain *L* is obtained. The concept of a number as developed herein comes laden with a seed of the necessary and absolutely infinite extension. Cantor uses *numerical amount, value*, and *limit* as equivalent.

Further, Cantor considers points on a line, defining the distance between them as a limit of a sequence and introducing relations of "greater than", "less than", and "equality". He introduces an axiom that, vice versa as well, a point on a line corresponds to each numerical value, the coordinate of such point being equal to this numerical value, and moreover, equal in the sense explained in this paragraph. Cantor calls this assertion as axiom, as it is not provable in its very nature. Thanks to this axiom, numerical values additionally gain definite objectivity on which, however, they do not depend a bit.

In accordance with the above, Cantor considers a point on a line as definite, if its distance from 0 considered with a definite sign is set as an λ-type numerical amount, value, or limit.

Further, Cantor defines multitude of points or point sets and introduces a concept of an accumulation point of the point set. A neighborhood shall be understood as any interval which contains this point. Thus, together with a set of points an ensemble of its accumulation points is defined. This set is known as the first derivative point set. If it consists of an infinite number of points, a second derivative point set may be formed of it, and so on [Cantor (1872)].

The introduction of the concept of an accumulation point (condensation point) was rewarding. Other mathematicians like H. Schwarz and U. Dini started using it right away. Cantor goes further in building a hierarchy of sets and displays the correspondence between sets. In the next five years, Cantor created his theory of point sets and proceeds to construction of a scale of infinities.

Definitions before Cantor were mostly constructive (how to construct an object) or functional (what for). In the humanities other definition types were used, i.e. negative and descriptive. Thanks to such definitions, one could start considering an object without knowing much about it, e.g., only one reference characteristic thereof. Further in the process of study, the volume of information grew. Any object could be separated from others with the help of a negative definition. Cantor adapted these principles of description to mathematical objects.

Objects in the theory of sets are constructed differently than before. Thanks to the high generality of construction, mathematics based on the theory of sets grew into an abstract science. Its notions did not depend on geometrical or physical meaning, or on applications.

Cantor wrote: "The process of correct formation of notions is, to my mind, the same everywhere: they take any item devoid characteristics, which is initially a pure and simple name or sign A, and add to it in a regular pattern various, even innumerable, straightforward predicates, the meaning of which is already known from ideas at hand and which shall not disagree. Thanks to this, one can determine relation of A to the already existing notions, especially to contiguous ones. Once this has been completed, all the conditions precedent to reviviscence of the dormant notion *A* are possessed, and this notion comes into being supplied with such intrasubjective reality as may be expected from notions in general. And then it is up to metaphysics to ascertain its transient notion". (Metaphysics was understood by Cantor as application areas) [Кантор (1985), c. 103–104].

Cantor analyses definitions provided by Weierstraβ, Dedekind and himself: "The definition of any surd real number is in correspondence with a strictly defined set of the first cardinality rational numbers. This is the common feature of all forms of definitions. The difference is in the time of creation with the help of which the set is joined with the number defined thereby and in the conditions to be met by the set to constitute an appropriate basis for respective definition of a number. In the first form [of Weierstraβ] the time of creation binding the set with the number defined thereby is in inception of sums. Dedekind uses the sum total of all rational numbers. However, it is accompanied by a large drawback – numbers in the analysis are never presented in the form of "sections" in which they have to be written using a pretty artificial and complicated way" [Кантор (1985), c. 81].

Introducing definitions of new notions, Cantor uses a different method of their formation for the first time: "In the theory set forth herein, a numeric value initially appears as something pointless, just like a constituent of theorems adding reality to it, for example, a theorem stating that respective array has this numeric value as a limit… The notion of the number as developed herein carries an embryo of the necessary and absolute generalization".

Weierstrass defined a real number as a limit of partial amounts of absolutely covergent series, noting the need in arithmetization of the concept of a limit. A point on a line corresponds to each number. However, it is not obvious that a number corresponds to each point. Cantor considers points on a line, defining the distance between them as a limit of a sequence and introducing relations of greater than, less than, and equal to. He introduces an axiom that, vice versa, a point on a line corresponds to each numerical value, the coordinate of such point being equal to this numerical value – equal in such meaning as set forth in this paragraph. Cantor calls this statement an axiom, as it is indefensible due to its very nature. Thanks to it, numerical values additionally gain definite objectivity on which, however, they do not depend at all. Dedekind believes that numbers are subjects of the "world of our thoughts", and it is our right to believe they are related to points. Unlike the above authors, he defines a real number as a limit of partial amounts of absolutely convergent series, noting the need in arithmetization of the concept of a limit.

Cantor was developing a perceptual theory of point sets, verily believing that applications were subsidiary issues. Years later, his theory of point sets devised as a summary of contemporary analysis formed basis of mathematics.

Dedekind developed and arithmetical concept of a number as an algebraist, not being inclined to problems of analysis. Fifteen years later, his design led to creation of Dedekind-Peano arithmetic axiomatics.

Heine pursued educational goals. His narrative on limit and continuity was included in modern courses of analysis. Simultaneously, he set forth a number of fundamental principles: on disregarding of a certain number of point, a covering lemma, concept of uniform convergence.

Creation of Charles Méray was recognized by his fellow-countrymen a century later and is now called a "Méray-Heine" or "Méray-Cantor concept of a number".

After Cantor created the set theory, the language and internal structure of mathematics changed. It did not anymore need geometrical or physical interpretation and gained a material descriptive component. Language and descriptive forms became the creating tool. The set theory was created as a continuation of arithmetics. However, already ten years later it formed basis of the theory of a real number. It provided the opportunity to analyze the finest shades of design of mathematical objects and links between them. Many definitions and statements were formed verbally, retaining a high abstraction degree. This caused to discussions among mathematicians devoted to contradictions many of which were linguistic in their nature. However, a new theory was created as a result, descriptive theory of multitudes, the key results of which belong to mathematicians of Warsaw and Moscow schools.

In the process of reasoning this notion becomes overgrown with features, fills with content, reserving a possibility of further development. Having introduced a definition of, say, enumerable sets, after some reasoning Cantor introduces a more accurate definition again. As opposed to Port-Royal logic, the volume and structure of a notion in the beginning of reasoning differs from that in the end. This principle derived from the humanities came to mathematics for the first time in Cantor's works. Thereafter, it formed the basis of the descriptive theory of sets. The principle of enhancement of a notion first defined concisely enabled development of mathematical ideas and raising problems in the interior of the very mathematics. Before the emergence of the theory of sets, problems used to arise out of practice. According to N. Lusin, the high generality of definitions and concepts of the theory of sets and differential calculus has allowed to "obtain many delicate facts which were not expected and were even impossible from the former point of view. Thus, assuming various sets of points to be an E set, we study a derivative on various sets, starting with continuous functions f (x) and, comparing numerical results at different E, we get into the deepest structural properties of continuous functions. These properties are kind of a microscopic study of behavior of the continuous function in the point x0 which disclose extraordinary richness of various correlations of new data which ultimately deeply affect the flow of function f (x) itself throughout the interval [a, b] ".

Thanks to the high generality of definitions in the theory of sets, there appeared a possibility to introduce a wide class of functions and geometrical forms in the study. Their description would have been impossible on the basis of the mathematical analysis of Cauchy and Weierstraβ.

The first period of development of the theory of sets was called the "naive" theory. Notions started to form verbally, new expressions derived from the already known ones by means of operations with words (terms) in accordance with grammar rules. At the same time the logic structure of the language was at variance with its grammatical structure, which caused paradoxes of the theory. Russell analyzed this process in his theory of descriptions distinguishing two types of relations of signs and the designated object – names and descriptions.

The descriptive theory of sets was conceived in works of Baire, Borel and Lebesgue in connection with the problem of measurability of sets. It was intended to study the internal structure of sets depending on operations by means of which these sets could be constructed from more simple, for example, closed or open subsets of Euclidean, metric or topological space. These operations include joining, meeting, taking the complement, etc.

N. Luzin, M. Suslin and P. Alexandrov in Russia, W. Sierpinski in Poland, and F. Hausdorff in Germany continued developing this theory.

The main principle of the descriptive theory of sets was laid down in 1915 by Luzin in his thesis: "Given: a structural property of a function. Find its analytic expression."

Great attention was given to classification of theorems depending on their dependence on the axiom of choice by W. Sierpinski in his works. His cooperation with N. Luzin was fruitful for both of them. Luzin wrote : "In 1916–1918 aspiring to fulfill the Lebesgue offered program of study of the most general sets, Suslin, W. Sierpinski and myself, who, so to say, came to study a new class of the point sets that is trivially goes beyond the class of sets and measurable B, however, originating from sets which can be defined without any transfinite numbers. Due to the close interrelation between these sets and series of polynomials, they were called analytical sets as offered by Lebesgue."

The issue on measurability of a set is important for the appendix to the analysis and the theory of functions. All known evidence of measurability of certain sets is connected with its structural, i.e. descriptive, properties.

Contemporary development of the descriptive theory of sets in common topological spaces is in demand in other areas of mathematics, e.g., in the theory of potential.

In 1934, in his article "Differential calculus", Luzin distinguishes the mathematical analysis of the high and low styles intended for different purposes – strict educational and creative research. It is possible to use the petitio principii in a research work. To the contrary, in courses of low style "excessive" notions are introduced, making important problems unsolvable due to the introduced notions. However, these problems are perfectly expressed in the language of this small analysis .

So a logical rigor of definitions is a good serve but a bad master.

In last third of the 19th – first third of the 20th century mathematics shaped into a new form – it gained the role of a fundamentally self-sufficient theory. The theory of sets, the theory of measure has played a significant role in this process. Remember, the degree of generality of many theorems varied to negligible small sets with measure zero. For example, in Russian indeterminacy is not a synonym of doubtfulness, it is more tolerant to contradictions. This gave the notions a chance to develop and come up with new, more exact definitions (for example, analytical and projective sets) in the course of research. This process was largely associated with the increase of the role of language as generating structure using grammatical means and philosophical sense. Thanks to this process the Moscow school of the theory of sets led by Luzin occupied leading positions.